\documentclass[12pt,a4paper]{elsarticle}
\usepackage{amssymb,amscd,amsmath,xypic}

\newtheorem{theorem}{Theorem}[section]
\newtheorem{lemma}[theorem]{Lemma}
\newtheorem{proposition}[theorem]{Proposition}

\newdefinition{definition}[theorem]{Definition}
\newdefinition{remark}[theorem]{Remark}
\newdefinition{remarks}[theorem]{Remarks}
\newdefinition{example}[theorem]{Example}

\newproof{pf}{Proof}

\newcommand\cO{{\mathcal O}}

\newcommand\cV{{\mathcal V}}

\newcommand\aq{/_{_{\operatorname{aff}}}\,}
\newcommand\gq{/\hspace*{-2pt}/}

\newcommand\End{\operatorname{End}}

\newcommand\Spec{\operatorname{Spec}}

\begin{document}
\title{Observable subgroups of algebraic monoids}

\author[ren]{Lex Renner\fnref{fn1}}
\ead{lex@uwo.ca}
\author[rit]{Alvaro Rittatore\fnref{fn2}}
\ead{alvaro@cmat.edu.uy}

\fntext[fn1]{Partially supported by a grant from NSERC.}
\fntext[fn2]{Partially supported by grants from IMU/CDE,
NSERC and PDT/54-02 research project.}

\address[ren]{University of Western Ontario, London, N6A 5B7,  Canada.}

\address[rit]{Facultad de Ciencias, Universidad de la Rep\'ublica,
  Igu\'a 4225,  11400 Montevideo, Uruguay.}

\begin{abstract}
A closed subgroup $H$ of the affine, algebraic group $G$ is called
\emph{observable} if
$G/H$ is a quasi-affine algebraic variety. In this paper we define the
notion of an
observable subgroup of the affine, algebraic \emph{monoid} $M$. We prove that a
subgroup $H$ of $G$ is observable in $M$ if and only if $H$ is closed
in $M$ and there are  ``enough'' $H$-semiinvariant functions in $\Bbbk
[M]$. We show also
that a closed, normal subgroup
$H$ of $G$ (the unit group of $M$) is observable in $M$ if and only if
it is closed in $M$.
In such a case there exists a \emph{determinant} $\chi : M \to \Bbbk$
such that 
$H\subseteq \operatorname{ker}(\chi)$. As an application, we show that
in this case the 
\emph{affinized quotient} $M\aq H$ of $M$ by $H$ is an affine
algebraic monoid 
scheme with unit group $G/H$.
\end{abstract}

\maketitle

\section{Introduction}

A closed subgroup $H$ of the affine algebraic group $G$ is called an
{\em observable subgroup} if the homogeneous space $G/H$ is a
quasi-affine variety. Such subgroups have been researched extensively,
notably by F.~Grosshans, see \cite{kn:grosshanslocal} for a survey on
this topic, and Theorem \ref{thm:obsercond} below for
other useful characterizations of observable subgroups. In \cite{kn:oaoag}
the authors presented the notion of an \emph{observable action} of $G$
on the affine variety $X$, together to its basic properties. In this
paper we develop further the notion of observable actions.  In
particular, we  investigate the situation where $M$ is an affine
\emph{algebraic monoid} with unit group $G$, and $H$ is a closed
subgroup of $G$, such that the action of $H$ on $M$ by left
multiplication is observable.   In this case, we say that
$H$ is \emph{observable in $M$}.

We describe now the organization of this paper. 

In Section \ref{sec2} we
provide the basic definitions and properties of observable actions and
affinized quotients. In Section \ref{sec3} we give several
characterizations of observable subgroups. 
In Theorem \ref{thm:istabcond} we deduce a number of useful
consequences from the assumption
that $H$ is observable in $M$. In particular it follows that $H$ is an
observable subgroup of $G$, and that it is closed in $M$.
In Theorem \ref{theo:1} we characterize the observable subgroups of $M$ in
terms of semiinvariants. In Theorem \ref{thm:obserrep} we show that
$H$ is an observable subgroup of $M$ if and only if $H$ is the
isotropy group of some vector $v\in V$ in some rational representation 
$\rho : M\to \End(V)$ of $M$. 
In the final section (Section \ref{sec4}) we use the results of the
previous section to study the affinized quotient of an affine 
algebraic monoid by a closed normal subgroup. 
In Theorem \ref{thm:normalimplisobse} we show that if $H$ is
a closed, normal subgroup of $G$, closed in $M$, then $H$ is observable in $M$.
Whether this is true for nonnormal closed subgroups of $G$ is an open question.
See Remark \ref{rem:openquest}. As an application, we show that the 
\emph{affinized quotient} of an affine algebraic 
monoid $M$ by a normal subgroup $H$, closed in $M$, is an algebraic
monoid, with unit group $G/H$.

\medskip

{\sc Acknowledgements: }
This paper was written during a stay of the second author at the
University of Western Ontario. He would like to thank them for the kind
hospitality he received during his stay.

\section{Preliminaries}
\label{sec2}

Let $\Bbbk $ be an algebraically closed field. We work
with affine algebraic varieties $X$ over $\Bbbk $. An algebraic
group is assumed to be a smooth, affine, group scheme of finite type
over $\Bbbk $. If $X$ is an affine variety over $\Bbbk $ we
denote by $\Bbbk [X]$ the ring of regular functions on $X$. If
$I\subset \Bbbk[X]$ is an ideal, we denote by $\mathcal
V(I)=\bigl\{x\in X\mathrel{:} f(x)=0 \ \forall \, f\in I\bigr\}$. If
$Y\subset X$ is a subset, we denote by $\mathcal I(Z)=\bigr\{ f\in
\Bbbk[X]\mathrel{:} f(y)=0\ \forall\, y\in Y\bigr\}$.  Morphisms
$\varphi:X\to Y$ between affine varieties correspond to morphisms of
algebras $\Bbbk[Y] \to \Bbbk[X]$, by $\varphi\mapsto \varphi^*$,
$\varphi^*(f)=f\circ \varphi$.  If $X$
is irreducible we denote by $\Bbbk (X)$ the field of rational
functions on $X$. If $A$ is any integral domain we denote by $[A]$ its
quotient field. Thus if $X$ is an irreducible affine variety, then
$\Bbbk(X)=\bigl[\Bbbk [X]\bigr]$.

Let $G$ be an affine algebraic group and let $X$ be an algebraic variety.
A \emph{(regular) action} of $G$ on $X$
is a morphism $\varphi:G\times X\to X$, denoted by $\varphi(g,x)=g\cdot
x$, such that $(ab)\cdot x= a\cdot (b\cdot x)$ and $1\cdot x=x $ for
all $a,b\in G$ and $x\in X$. Since all the actions we  work with
are regular, we  drop the adjective regular. The \emph{orbit} of $x\in X$ is
denoted by $\cO (x)=\{ g\cdot x\mathrel{:} g\in G\}$.

If $G\times X\to X$ is a regular left action of $G$ on $X$ we consider the
induced right   action of $G$ on $\Bbbk [X]$, defined as
follows. If $f\in
\Bbbk [X]$ and $g\in G$,  then
 $(f\cdot g)(x)=f(gx)$. It is
well known that $G$-stable
closed subset of $X$ correspond to $G$-stable radical ideals  of
$\Bbbk[X]$. We say that
$f\in\Bbbk [X]$ is \emph{G-invariant} if $f\cdot g=f$ for any $g\in
G$.  The set of
$G$-invariants ${}^G\Bbbk [X]$ forms a
$\Bbbk $-subalgebra of $\Bbbk [X]$, possibly non-finitely
generated. Analogous considerations can be made if we start with a
right action $X\times G\to X$.

A \emph{finite dimensional (rational) $G$-module} is a
finite dimension $\Bbbk$-vector space $V$ together with a left action
of $G$ on $V$ by linear automorphisms. A right action of
$G$ on $V$ defines, in a similar way, the
notion of a \emph{right $G$-module}.

Recall that an \emph{algebraic monoid} $M$ is an algebraic variety
together with an associative product $m:M\times M\to M$ with neutral
element $1$, such that $m$ is a morphism of algebraic varieties.  We
denote the \emph{set of idempotent elements} of $M$ by
$E(M)=\bigl\{e\in M\mathrel{:} e^2=e\bigr\}$. We denote the \emph{unit
  group} of $M$ by $G(M)=\bigl\{ g\in M\mathrel{:} \exists \,
g^{-1}\in M\,,\ g^{-1}g=gg^{-1}=1\bigr\}$. It is known that
$G(M)$ is an algebraic group, open in $M$ (see \cite{Re05} and
\cite{Ri98}).

\subsection{Characters of affine algebraic monoids}

In this section we establish the basic facts about {\em extendible
characters} for the case of linear algebraic monoids.

\begin{definition}
Let $M$ be an algebraic monoid. A {\em character}\/  of
$M$ is a  morphism of algebraic monoids  $M\rightarrow \Bbbk$.
 We denote the monoid of characters of $M$ by
\[
\mathcal X(M)=\bigl\{ \chi\in \Bbbk[M]\mathrel{:}
\chi(ab)=\chi(a)\chi(b)\ \forall\, a,b\in M\,,\ \chi(1)=1\bigr\}.
\]
If $G=G(M)$, then restriction induces an injective morphism of (abstract)
monoids $\mathcal X(M)\hookrightarrow \mathcal X(G)$.
\end{definition}

If $M$ is an irreducible affine algebraic monoid, then there exists
$n\geq 0$ and a morphism of algebraic monoids $\rho:  M\hookrightarrow
\operatorname{M}_n(\Bbbk)$, such that $\rho$ is closed
immersion (see for example \cite[Theorem 3.8]{Re05}). This motivates
the following
definition.

\begin{definition}
Let $M$ be an irreducible affine algebraic monoid. A character $\chi
\in \mathcal X(M)$ is called a \emph{determinant} if
$\chi^{-1}(0)=M\setminus G(M)$.

By the considerations above, determinants always exists.
\end{definition}

\begin{remark}
(1) Observe that if $\det M\to \Bbbk$ is a determinant, then
$G(M)=M_{\det}$. In particular $\Bbbk\bigl[G(M)\bigr]=\Bbbk[M]_{\det}$.

\noindent (2) Clearly, $\mathcal
X(M)\subset
\mathcal X \bigl(G(M)\bigr)$
is a unital submonoid that generates $\mathcal X\bigl(G(M)\bigr)$ as
a group. Indeed, let $\det :M\to \Bbbk $ be a determinant and $\chi \in
\mathcal X\bigl(G(M)\bigr)$. Then there exists $f\in \Bbbk[M]$ and
$n\geq 0$ such that $\chi=\frac{f}{\det^n}$. It follows that $f|_{_G}$
is a character, and hence, by continuity, $f\in \mathcal X(M)$.
\end{remark}

The notion of an {\em extendible
character} is very useful in the study of observable subgroups of algebraic
groups. We extend this notion to the setting of algebraic monoids.

\begin{definition}
Let $M$ be an affine algebraic monoid and let $H\subset G(M)$ be a
closed subgroup. A non-trivial  character $\chi\in \mathcal X(H)$
is \emph{extendible} if there exists a non-zero \emph{semiinvariant
  element of
weight $\chi$}; that is,  there exists
$f\in k[M]$ such that
$f\cdot x=\chi(x)f$ for all $x\in H$.  Such an element is called an
{\em extension of $\chi$}. We denote the monoid of extendible
characters of $H$ by $E_M(H)$.

Clearly, restriction to $G=G(M)$ induces an injective
homomorphisms of (abstract) monoids $E_M(H)\hookrightarrow E_G(H)$.
\end{definition}

\begin{remark}
Observe that if $\chi$ is an extendible character and $f$ is an
extension of $\chi$, we can suppose that $f(1)=1$. Then $f(x)=(f\cdot
x)=(1)\chi(x)f(1)=\chi(x)$ for all $x\in N$, and thus $f$ is an
extension of $\chi$ to a regular function of $M$.
\end{remark}

\begin{definition}
Let $M$ be an affine algebraic monoid and $V$ a finite rational
$M$-module. For every $\alpha \in
V^{*}$ and $v \in V$ we define  $\alpha | v :M \rightarrow \Bbbk $ as
$(\alpha|v)(x)=\alpha(x\cdot
v)$ for all $x \in M$. A function of this form is called an
{\em $V$--representative function}\/  or simply a
{\em representative function.}
\end{definition}

\begin{definition}
Let $G$ be an algebraic group and let  $\chi\in \mathcal X(G)$ be a
character. Let $V$ be 
a right $G$-module,  with action $\varphi:V\times G\to V$,
$\varphi(v,g)=v\cdot g$. Then the \emph{twisted representation}
$V_\chi$ is 
defined as follows. As a vector space $V=V_\chi$, and the action
$V_\chi\times G\to V_\chi$ is given by $v\star
g=\chi(g)(v\cdot g)$, for every $g\in G$, $v\in V_\chi$.
\end{definition} 

The following theorem is a straightforward generalization of the corresponding
(well known) result for algebraic groups (see \cite[Theorem 7.2.3]{fer-ritt}).
We include a proof for the sake of completeness.

\begin{theorem}
\label{theo:h1}
 Let $M$ be an affine algebraic monoid with unit group $G$,
 $H\subset G$ a closed subgroup
 and  $V$  a finite dimensional rational right $H$-module. There
 exists a
 finite dimensional rational right $M$-module $W$, an extendible character
 $\chi: H \to G_m$ and an injective morphism $\iota : V \to
(W|_{H})_{\chi^{-1}}$.
\end{theorem}

\begin{pf} 
Given $V$ as above, we proceed as in the proof of \cite[Theorem
7.2.3]{fer-ritt}.  It is well known that there exists an injective
morphism of $H$-modules
$\theta : V \to \bigoplus_{I} \Bbbk [H]$, where $I$ is a finite set of
indexes (see for example \cite[Theorem 4.3.13]{fer-ritt}).
Consider the $H$-morphism $\alpha=\bigoplus \pi : \bigoplus_{I} \Bbbk [M] \to
\bigoplus_{I} \Bbbk [H]$, where   $\pi : \Bbbk [M] \to \Bbbk [H]$ is
the canonical
projection.

Call $V'=\alpha^{-1}\bigl(\theta(V)\bigr)\subset \bigoplus_{I} \Bbbk [M]$ and
$\beta$ the restriction of $\alpha$ to the $H$--submodule $V'$:
\[
\xymatrix
{
V' \ar@{->>}[d]_-{\beta} \ar@{^{(}->}[r] &
\bigoplus_{I} \Bbbk [M] \ar@{->>}[d]^-{\alpha} \\
V \ar[r]_-{\theta} & \bigoplus_{I} \Bbbk [H]
}
\]

Let ${\mathcal F}$ be a finite $\Bbbk $-linear basis of $V$ and let
${\mathcal F}_{0} \subset V'$ be a finite set such that
$\beta(\mathcal F_{0})={\mathcal F}$. Call
$R$ the finite dimensional $G$-submodule of $\bigoplus_{I} \Bbbk [M]$ generated
by ${\mathcal F}_{0}$; then $R$ is an $M$-module. Let $S$ be finite
dimensional  $H$-submodule
of $R|_{_H}$ (contained in  $V'$) generated by ${\mathcal
  F}_{0}$. In this way, we produce an exact sequence of rational
$H$--modules $0 \to U \to S \to V \to 0$. Call $n ={\operatorname
  {dim}}_{\Bbbk } U$ and consider the commutative diagram below
 \[
\def\objectstyle{\textstyle}
\def\labelstyle{\scriptstyle}
\xymatrix{
0 \ar[r] & U \otimes \bigwedge^nU \ar[r] \ar[d] & S \otimes
\bigwedge^nU \ar[r] \ar[d] & V \otimes \bigwedge^nU
\ar@{.>}[dl]^-{\varphi} \ar[r] & 0 \\
0 \ar[r] & U \land \bigwedge^nU=0 \ar[r] &  S \land \bigwedge^nU  & &
}
\def\objectstyle{\textstyle}
\def\labelstyle{\textstyle}
\]

\noindent where all the solid arrows are the canonical ones,
the first row is exact and in the second row
the term $U \land \bigwedge^nU$ equals zero by dimensional reasons
--- notice that all the exterior products are
taken inside the exterior algebra $\bigwedge R$.
One can now prove that $\varphi $ is bijective (see the proof of \cite[Theorem
7.2.3]{fer-ritt}).

As $S$ and $U$ are $H$-submodules
of $R|_{_H}$, we can view $\varphi$ as an injective $H$--morphism
$\varphi_{1} :  V \otimes \bigwedge^{n}U \to  \bigwedge^{n+1}R|_{_H}$.
Call $\chi$ the rational character
associated to the one dimensional $H$-module $\bigwedge^{n}U$, i.e.\ the
character defined by the formula $x\cdot u = \chi(x) u$ for all $x \in H$.
The map $\varphi_{2} :  V \to  \bigwedge^{n+1}R$,
$\varphi_{2}(m)=\varphi_{1}(m \otimes u)$, satisfies that for all $x \in H$,
\[
\begin{split}
\varphi_{2}(m\cdot x) &= \varphi_{1}\bigl((m\cdot x) \otimes u\bigr) =
\varphi_{1}\bigl(\chi^{-1}(x) (m\cdot x \otimes  u\cdot x)\bigr) \\
& =\chi^{-1}(x)\varphi_{1}\bigl( (m \otimes  u)\cdot x\bigr) =
\chi^{-1}(x) \varphi_{1}(m \otimes  u)\cdot x\\
& =
\chi^{-1}(x) \varphi_{2}(m)\cdot x.
\end{split}
\]

Hence, if we let $W= \bigwedge^{n+1}R$ and $\iota= \varphi_{2}$, the
proof of the theorem will be complete once we prove that the character
$\chi$ is extendible.
Consider the $M$--module $\bigwedge^{n}R$ and let $\alpha \in
(\bigwedge^{n}V)^{*}$ be such that $\alpha(u)=1$.
 Then the representative function $\alpha |u \in \Bbbk [M]$ is an extension of
 $\chi$.  Indeed,  if $x \in H$ and $m\in M$, then
\[
\bigl((\alpha  |u)\cdot x\bigr)(m) = (\alpha |u)(xm)= \alpha
\bigl(u\cdot (xm)\bigr) = \chi(x) \alpha (u\cdot m
)=\chi(x)(\alpha|u)(m)
\]
and $(\alpha
|u)(1)= \alpha(u) =1$.  So that $\alpha |u$ is an extension of
$\chi$.
 \qed \end{pf}

Let $M$ be an affine algebraic monoid with unit group $G$. 
Many results about extendible characters for subgroups of $G$
apply to the situation of closed subgroups of $M$.

\begin{lemma}\label{lema:1}
Let $M$ be an affine algebraic monoid with unit group $G$,  and
$H\subset G$ a subgroup, closed in $M$.

\noindent (1) Let  $\rho\in E_M(H)$  be an extendible character and let $f \in
k[M]$  be a semiinvariant of weight $\rho$.  Consider the left action
$M\times \Bbbk[M]\to \Bbbk[M]$ given by $(x\cdot f)(m)=f(mx)$, $x,m\in
M$. Then for any $x \in M$, $x\cdot f$ is also a semiinvariant of weight
$\rho$.

\noindent (2) If $\pi: \Bbbk[M] \rightarrow  \Bbbk[H]$ is the canonical
projection and
$\rho \in E_M(H)$ is an extendible character, there
exists an extension $f \in
k[M]$ of $\rho$ such that $\pi(f)=\rho$.

\noindent (3) A character $\rho\in \chi(H)$ is extendible if and
only if there exists a
rational right $M$--module $V$ and an injective morphism of right $H$--modules
$\iota : \Bbbk\rho \rightarrow V$. In other words $\rho$ is
extendible if and only if there exists a rational right $M$--module $V$ and
a non-zero element $v \in V$ such that $v\cdot x = \rho(x)v$ for all
$v \in H$. Moreover, the right $M$--module $V$ can be taken
to be a finite dimensional right $M$--submodule of $\Bbbk[M]$.

\noindent (4) For any $\gamma \in {\mathcal X}(H)$ there exists $\rho \in
E_{M}(H)$ such that $\gamma \rho \in E_{M}(H)$.
\end{lemma}
\begin{pf} 
Many of these results are straightforward generalizations of the 
corresponding results for algebraic groups. See \cite[Lemma 7.2.8]{fer-ritt}.
For the sake of completeness we include the proofs of (3) and (4).

\noindent To prove (3) suppose that $\rho$ is extendible, let $f \in \Bbbk [M]$ 
be an extension and call $W$ the rational right $G$--submodule of $\Bbbk
[M]$ generated
by $f$. Then the map
$\iota:\Bbbk \rho \to W$, $\iota(\rho)=f$, does the job.
Conversely, if one has an injective morphism $\iota:\Bbbk \rho \to
W|_{_H}$  and call $n = \iota(\rho)$, then $n\cdot x=
\iota(\rho)\cdot x= \iota( \rho\cdot x)
= \rho(x) \iota(\rho)=\rho(x)n$. Take now $\alpha \in N^{*}$ such that
$\alpha(n)=1$ and consider $f= \alpha | n$. Then, $f\cdot x(m) = (\alpha|n)(xm)
\alpha  \bigl(n\cdot (xm)\bigr)=\rho(x) (\alpha| n)(m) =\rho(x)f$ for
all $x \in H$, and $f(1)=
(\alpha |n) (1)= \alpha(n)=1$.

\noindent To prove (4) consider $V=\Bbbk  \gamma$. By Theorem
\ref{theo:h1} there exists an extendible
character $\rho$, a right $M$--module $W$ and an injective $H$--morphism
$\iota: \Bbbk \gamma \to \bigl(W|_{_H}\bigr)_{\rho^{-1}}$. Call
$n=\iota(\gamma)$,
then for any $x \in H$ we have that $n\cdot x =  \iota(\gamma)\cdot x=
\rho(x)\iota( \gamma\cdot x)= \rho(x)\gamma(x)\iota(\gamma)=
\rho(x)\gamma(x)n$.
From (3) we conclude that $\rho\gamma \in E_{M}(H)$.
 \qed \end{pf}

\subsection{Observable actions of algebraic groups}

Observable subgroups were introduced by
Bialynicki-Birula, Hochs\-child and Mostow in \cite{kn:obsdef}. Since then
they have been researched extensively,
notably by F.~Grosshans. See \cite{kn:grosshanslocal} for a survey on
this topic. This basic notion has recently been generalized and reformulated in the 
context of geometric invariant theory  by the authors. See \cite{kn:oaoag}.

\begin{definition}
\label{defi:obsgrp}
Let $G$ be an affine algebraic group and let $H\subset G$ be a closed subgroup. The
subgroup $H$ is {\em observable} in {$G$}\/ if and only if for any
nonzero $H$-stable ideal $I \subset \Bbbk [G]$ we have that $I
^G\neq (0)$, where we consider the right action of $H$ on $\Bbbk[G]$ given
by $ (f\cdot h)(a)=f(ha)$, for all $f\in \Bbbk[G]$ $h\in H$ and $a\in G$.

Analogously, one can define the notion of a \emph{right
observable subgroup}, by considering the action $H\times
\Bbbk[G]\to\Bbbk[G]$,  $(h\cdot f)(a)=f(ah)$ for $f\in \Bbbk[G]$,
$a\in G$, $h\in H$.
\end{definition}

\begin{example}
(1) If $U\subset G$ is a closed unipotent subgroup, then $U$ is observable,
since any $U$-module has non-zero invariant elements.

\noindent (2) If $H\subset G$ is a normal closed subgroup, then $H$ is
observable.  This follows from condition (2) of Theorem
\ref{thm:obsercond} below.

\noindent (3) Let $H\subset G$ be a closed subgroup, such that
$\mathcal X(H)=1$. Then $H$ is observable in $G$. This follows for
example from condition (6) of Theorem
\ref{thm:obsercond} below.
\end{example}

We now present a collection of equivalent definitions of observability. The proofs
can be found in
\cite[Thms.~10.2.9 and 10.5.5]{fer-ritt}. We give a different proof
for the fact that $H$ is observable in $G$ if and only if $E_G(H)$ is
a group (equivalence  (1)
$\Longleftrightarrow$ (7) of Theorem \ref{thm:obsercond}). This will
provide some insight into the more general situation of algebraic monoids given in
Theorem \ref{theo:1}.

\begin{theorem}
\label{thm:obsercond}
Let $G$ be an affine algebraic
group and $H\subset G$ a closed subgroup. Then the following
conditions are equivalent:

\begin{enumerate}
\item The subgroup $H$ is observable in $G$.

\item The homogeneous space $G/H$ is a quasi--affine variety.

\item For an arbitrary proper and closed subset $C \subsetneq
G/H$, there exists an non-zero invariant regular function $ f \in
\Bbbk [G]^H$ such that
$f(C)=0$.
\end{enumerate}

\noindent Moreover, if  $G$ is connected the above conditions 
are equivalent to any of the following.

\begin{enumerate}
\item[(4)]  $H= \bigl\{x \in G\mathrel{:} x\cdot f = f, \ \forall f \in
\Bbbk [G]^H\bigr\}$.

\item[(5)] $\bigl[ \Bbbk [G]\bigr]^H=
\bigl[\Bbbk [G]^H\bigr]$.

\item[(6)] $E_G(H)=\mathcal X(H)$.

\item[(7)] $E_G(H)$ is a group.
\end{enumerate}
\end{theorem}

\begin{pf} 
To prove that (1)
$\Longrightarrow$ (7), let  $H$ be observable in $G$ and $\chi\in
E_G(H)$ an extendible character. Let $f\in \Bbbk[G]$ be a
semiinvariant of weight $\chi$. Then the ideal $I=f \Bbbk[G]$ is
$H$-stable, and hence there exists $g\in \Bbbk[G]$ such that $
(fg)\cdot x=fg$ for all $x\in H$. It follows that $g\cdot x =
\chi^{-1}(x)g$ for all $x\in H$. Thus $\chi^{-1}\in E_G(H)$.

To prove that (7)
$\Longrightarrow$ (1), let $I\subset \Bbbk[G]$ be a non-zero $H$-stable
ideal, and consider $H_{\operatorname{uni}}=R_u(H)[H,H]$. Then
$H_{\operatorname{uni}}$ is normal in $H$, with $\mathcal
X(H_{\operatorname{uni}})=\{1\}$, and such that
$H_{\operatorname{uni}}\backslash H$ is a torus; in particular $\mathcal
X(H)=\mathcal X(H_{\operatorname{uni}}\backslash H)$, since $\mathcal
X(H_{\operatorname{uni}})=\{1\}$.
It follows that $H_{\operatorname{uni}}$ is observable in $G$, and thus
$I^{H_{\operatorname{uni}}}\neq (0)$ is a $H$-module. Then
$I^{H_{\operatorname{uni}}}$ is a right
$(H_{\operatorname{uni}}\backslash H)$-module,
and hence
\[
(0)\neq I^{H_{\operatorname{uni}}}= \bigoplus_{\chi\in \mathcal
  X(H_{\operatorname{uni}}\backslash H)}I^{H_{\operatorname{uni}}}_\chi=
\bigoplus_{\chi\in \mathcal
  X(H)}I^{H_{\operatorname{uni}}}_\chi
\]
Hence, there exists $\chi \in \mathcal X(H)$ and $f\in I\setminus
\{0\}$ such that
$f\cdot x=\chi(x) f$ for all $x\in H$. Thus $\chi\in E_G(H)$. Since
$E_G(H)$ is a
group, then $\chi^{-1}\in E_G(H)$ and hence there  exists $g\in
\Bbbk[G]$ such that $g\cdot x=\chi^{-1}(x) g$ for all $x\in H$. It
follows that $0\neq fg\in I^H$.
 \qed \end{pf}

\begin{definition}  \label{obser.def}
Let $G$ be an affine algebraic group acting on an affine variety
$X$. We say that  the action is \emph{observable}
if for any nonzero $G$-stable ideal $I\subset
\Bbbk[X] $ we have that $I^G\neq (0)$.
\end{definition}

\begin{example}
(1) The action of an unipotent group on an affine variety is always
observable, since for every right $U$-module $M$, we have that $M^U\neq 0 $.

\noindent (2) Let $G$ be an algebraic group and $H$ a closed subgroup. Then
the action of $H$ on $G$ by left translations is observable if and only if
$H$ is observable in $G$, in the sense of Definition \ref{defi:obsgrp}, if and only
if the action by right translations is observable.
\end{example}

Observable actions have been studied in detail in
\cite{kn:oaoag}, where
several characterizations of this important property were presented. We
collect here the ones that we need in what follows. We begin by recalling
some notation.

\begin{definition}
Let $G$ be an affine group acting on an affine variety $X$. We let
\[
\Omega(X)=\bigl\{ x\in X\mathrel{:} \dim \cO(x) \text{ is maximal
and } \overline{\cO(x)}=\cO(x)\bigr\}.
\]
That is, $\Omega(X)$ is the set of orbits of maximal dimension that
are closed. The reader should be aware that $\Omega(X)$ can be empty.
\end{definition}

\begin{theorem}
\label{thm:obserchar}
Let $G$ be a connected  affine algebraic group acting on an
irreducible affine variety $X$. Then the following are equivalent

\begin{enumerate}
\item The action is observable.
\item
  \begin{enumerate}
  \item[(i)] $\bigl[\Bbbk[X]^G\bigr]=\bigl[\Bbbk[X]^G\bigr]$ and
  \item[(ii)] $\Omega(X)$ contains a non-empty open subset.
  \end{enumerate}
\item There exists  a nonzero invariant $f\in \Bbbk[X]^G$ such that  the action of
      $G$ on $X_f$ is   observable.
\end{enumerate}
\end{theorem}

\begin{pf} 
See \cite[Proposition 3.2 and Theorem 3.10]{kn:oaoag}. 
 \qed \end{pf} 









If $G$ is a reductive group  acting on an affine variety $X$, then a
result of Popov (see \cite[Theorem 
4]{Pop}) guarantees that   if $\Omega(X)$
is nonempty 
then it is an open set. In \cite{kn:oaoag} this result was used to
prove that in this case, 
condition (2ii) of Theorem \ref{thm:obserchar} is sufficient to guarantee
observability.

\begin{proposition}
\label{prop:obsred}
Let $G$ be a reductive group acting on an affine variety $X$. Then the
action is observable if and only if $\Omega(X)\neq \emptyset$.
\end{proposition}
\begin{pf} 
See    \cite[Theorem 4.7]{kn:oaoag}.
 \qed \end{pf}

\subsection{The Affinized quotient}

Let $G$ be an affine algebraic group acting on the affine variety $X$.
It is well known that the categorical quotient does not necessarily
exist, even when $\Bbbk[X]^G$ is finitely generated. However, if
$\Bbbk[X]^G$ is finitely generated, then
$\Spec\bigl(\Bbbk[X]^G\bigr)$ satisfies a universal property in
the category on the affine algebraic varieties.

\begin{definition}
Let $G$ be an affine algebraic group acting on an affine variety
$X$, in such a way that $\Bbbk[X]^G$ is finitely generated. The
\emph{affinized quotient} of the action is the morphism $\pi:X\to
X\aq G=\Spec\bigl(\Bbbk[X]^G\bigr)$.

It is clear that $\pi$ satisfies the following universal property.

{\em Let $Z$ be an affine variety and $f:X\to Z$  a morphism constant
  on the $G$-orbits. Then
  there exists a unique $\widetilde{f}:X\aq G\to Z$ such that the
  following diagram is commutative.}
\begin{center}
\mbox{
\xymatrix{
X\ar@{->}[r]^f\ar@{->}[d]&Z\\
X\aq G\ar@{->}[ur]_{\widetilde{f}}&
}
}
\end{center}
\noindent Indeed, it is clear that the induced morphism $f^*:\Bbbk[Z]\to
\Bbbk[X]$ factors trough $\Bbbk[X]^G$.
\end{definition}

\begin{remark}
\label{rem:quotnotsurj}
It is clear that the  morphism $\pi: X\to X\aq G$ is dominant.
However, $\pi$ is not necessarily surjective. For example, consider a
semisimple group $X=H$ and its maximal unipotent subgroup $G$
acting on $X$ by left translation.
\end{remark}






\section{Observable subgroups of algebraic monoids}
\label{sec3}

We now adapt the notion of observability to the situation 
of subgroups of algebraic monoids.

\begin{definition}
\label{defi:obsmon}
Let $M$ be an algebraic monoid with unit group $G$, and let $H\subset G$
be a closed subgroup. We say that $H$ is \emph{(left) observable in $M$} if the
action by left multiplication $H\times M\to M $, $h\cdot m=hm$, is
observable in the sense of Definition \ref{obser.def}.

Similarly, we say that $H$ is \emph{right observable in $M$} if the
action by right multiplication $M\times H\to M $, $m\cdot h=mh$, is
observable.
\end{definition}

\begin{example}
\label{exam:obsmon}
(1) Since the action of a unipotent group on an affine variety is
observable  (all orbits are closed), it follows that any unipotent subgroup
$U$ of $G(M)$ is observable in $M$.

\noindent (2) If $M=G(M)$ is an algebraic group, then Definitions
\ref{defi:obsgrp} and \ref{defi:obsmon} coincide.

\noindent (3) If follows from Proposition \ref{prop:obsred} that if
$H\subset G(M)$
is reductive, and closed in $M$, then $H$ is observable in $M$.
Indeed, if $H=\overline{H}$, then $\Omega(M)\neq \emptyset$.
\end{example}

\begin{theorem}
\label{thm:istabcond}
Let $M$ be an affine irreducible algebraic monoid with unit group $G$,
and let $H\subset M$ be an
observable subgroup.  Then

\noindent (1) The subgroup $H$ is observable in $G$.

\noindent (2) $H$ closed in $M$.

\noindent (3) $\bigl[ \Bbbk [M]\bigr]^H=\bigl[ \Bbbk [G]\bigr]^H=
\bigl[\Bbbk [G]^H\bigr]=
\bigl[\Bbbk [M]^H\bigr]$.

\noindent (4) The subgroup $H$ satisfies
\[
\begin{split}
H= &\ \bigl\{x \in G\mathrel{:} f\cdot x = f, \ \forall f \in
\Bbbk [G]^H\bigr\}= \\
 &\ \bigl\{x \in G\mathrel{:} f\cdot x = f, \ \forall f \in
\Bbbk [M]^H\bigr\}=\\
& \ \bigl\{x \in M\mathrel{:} f\cdot x = f, \ \forall f \in
\Bbbk [M]^H\bigr\}^0.
\end{split}
\]
Recall that if $N$ is an algebraic monoid, then $N^0$ is the unique
irreducible component containing $1$, see for example \cite[Thm.~4]{Ri06}.
\end{theorem}

\begin{pf} 
\noindent (1)
Let $C\subset G$ be a $H$-stable closed subset. Then
$\overline{C}\subset M$ is $H$-stable, and it follows that there
exists a $f\in \cV(\overline{C})^H\setminus \{0\}\subset
\Bbbk[M]$. Since $G$ is open in $M$, it follows that $f\in
\cV(C)^H\setminus \{0\}\subset \Bbbk[G]$.

\noindent (2) Since $H$ is observable in $M$, it follows from Theorem
\ref{thm:obserchar}  that
$\Omega(M)$ contains a non-empty open subset. Since $G=\bigcup Hg$, $G$ is
contained in
$M_{\max}$, and hence there exists $g\in G$ such that $g\in
\Omega(M)\cap G$, i.e.~such that $Hg$ is closed
in $M$. Since multiplication by an element of $g$ is an isomorphism,
it follows that $H$ is closed in $M$.

\noindent  (3) First observe that $\bigl[\Bbbk [M]^H\bigr]\subset
\bigl[ \Bbbk
[M]\bigr]^H=\bigl[ \Bbbk [G]\bigr]^H$.   Let $g\in
\bigl[\Bbbk[M]\bigr]^H$, and consider the
ideal $I=\bigl\{ f\in \Bbbk[M]\mathrel{:} fg\in\Bbbk[M]\bigr\}$. Then
$I$ is a non-zero $H$-stable ideal, and hence there exists $h\in
I^H\setminus \{0\}$. It follows that $r=hg\in \Bbbk[M]^H$ and hence
$g=\frac{r}{h}\in \bigl[\Bbbk[M]^H\bigr]$.
The remaining equality follows from
 Theorem
\ref{thm:obsercond}.

\noindent (4) The first equality follows  from Theorem
\ref{thm:obsercond}.  It is clear that $
A=\bigl\{x \in G\mathrel{:} f\cdot x = f, \ \forall f \in
\Bbbk [G]^H\bigr\} \subset  B=\bigl\{x \in G\mathrel{:} f\cdot x = f,
\ \forall f \in
\Bbbk [M]^H\bigr\}$. Let $x\in B$ and $f\in \Bbbk[G]^H$. By (3), it
follows that $f\in
\bigl[\Bbbk[G]^H\bigr]=\bigl[\Bbbk[M]^H\bigr]$. Hence, there exist
$g,h\in \Bbbk[M]^H$ such that $f=\frac{g}{h}$. It follows that $f\cdot
x=\frac{g\cdot x}{h\cdot x}=\frac{g}{h}$.

In order to prove the last equality, we first observe that $N=
\mathcal V\bigl(
\bigl\{x\mapsto  (f\cdot x)(m)
-f(m)\,, \ m\in M\,,\ f \in
\Bbbk [M]^H\bigr\} \bigr)$ is a closed submonoid of $M$. Since $G$ is open
in $M$, it follows that
\[
H=\overline{H}= \overline{\overline{H}\cap G}=N^0.
\]
 \qed \end{pf}

\begin{theorem}
\label{theo:1}
Let  $M$ be an irreducible affine algebraic monoid with unit group $G$  and
let $H\subset G$ be a subgroup, closed in $M$. Then the following conditions
are equivalent.

\noindent (1) The subgroup $H$ is observable in $M$.

\noindent (2) $E_M(H)$ is a group.  That is,  for every $\rho \in E_{M}(H)$,
$\rho^{-1} \in E_{M}(H)$.

\noindent (3) $E_{M}(H)= \operatorname{\mathcal X}(H)$, i.e.\ every rational
character is extendible.




\noindent (4) For every finite dimensional rational right $H$-module $V$
there exists a finite dimensional rational right $M$-module $W$ and an
injective morphism of $H$-modules $\xi: V \to W|_{_H}$.



\noindent (5) $H$ is observable in $G$ and for any (some) determinant
$\det: M\to \Bbbk$, then $\frac{1}{\det}|_{_H}\in E_M(H)$.
\end{theorem}

\begin{pf} 
In order to prove that (1) implies (2), assume that $H$ is observable
in $M$, and let $\chi\in E_M(H)$ be an
extendible character. Let $g$ be an extension of $\chi$. Then the ideal
$g \Bbbk[M]$ is a nonzero
$H$-stable ideal and hence there exists $f\in \Bbbk[M]$ such that $
(gf)\cdot x=gf$ for all $x\in H$. It follows that $f\cdot x =
\chi^{-1}(x)f$ for all $x\in H$. That is, $\chi^{-1}\in E_M(H)$.

Since $E_M(H)$ generates $\mathcal X (H)$ as a group, it is clear that
conditions (2) and (3) are equivalent.

To prove that (2) implies (1), assume  that $E_M(H)$ is a
group and let  $I\subset \Bbbk[M]$ be a non-zero $H$-stable
ideal. Let $H_{\operatorname{uni}}=R_u(H)[H,H]$, as in the  proof of
Theorem \ref{thm:obsercond}. Then $H_{\operatorname{uni}}$ is
observable in $H$. Consider  a determinant $\det:M\to
\Bbbk$. Since  $\mathcal
X(H_{\operatorname{uni}})=1$, it follows that
$H_{\operatorname{uni}}\subset \det^{-1} (1)$.  Thus, we have found an
$H_{\operatorname{uni}}$-invariant regular function $\det\in
\Bbbk[M]^{H_{\operatorname{uni}}} $, such that
$\Bbbk[G]=\Bbbk[M]_{\det}$. It follows from Theorem
\ref{thm:obserchar} that $H_{\operatorname{uni}}$ is observable in
$M$, and thus the right $H$-module
$I^{H_{\operatorname{uni}}}$ is not trivial. Then, as in the
proof of Theorem \ref{thm:obsercond}, it follows that
$I^{H_{\operatorname{uni}}}$ is a right  $(H_{\operatorname{uni}}\backslash
H)$-module,
and hence
\[
(0)\neq I^{H_{\operatorname{uni}}}= \bigoplus_{\chi\in \mathcal
  X(H_{\operatorname{uni}} \backslash H)}I^{H_{\operatorname{uni}}}_\chi=
\bigoplus_{\chi\in \mathcal
  X(H)}I^{H_{\operatorname{uni}}}_\chi.
\]
Hence, there exists $\chi \in \mathcal X(H)$ and $f\in I\setminus
\{0\}$ such that
$f\cdot x=\chi(x) f$ for all $x\in H$. Thus $\chi\in E_M(H)$. Since
$E_M(H)$ is a
group, it follows that  $\chi^{-1}\in E_M(H)$. Thus, there  exists $g\in
\Bbbk[G]$ such that $g\cdot x=\chi^{-1}(x) g$ for all $x\in H$. It
follows that $0\neq fg\in I^H$.

To prove that (2) implies (4), we follow the idea of the
proof of \cite[Theorem 10.2.9]{fer-ritt}. If $V$ is a finite dimensional
rational right $H$-module, using
Theorem \ref{theo:h1} we deduce the existence of
an extendible
character $\rho$, a finite dimensional rational right $M$-module $W$ and
an injective map $\iota: V \to W$ such that $\iota(m\cdot x)
= \rho^{-1}(x)  \iota(m)\cdot x$ for all $x \in H$. By hypothesis the character
$\rho^{-1}$ is extendible, so if we take $f \in \Bbbk [M]$ that extends
$\rho^{-1}$ and call $W_{0}$ the right $M$-submodule of $\Bbbk [M]$ generated
by $f$, we can define an injective map $\xi: V \to W \otimes
W_{0}$,
given as $\xi(m)=\iota(m) \otimes f$. If we endow $W \otimes W_{0}$
with the diagonal right $M$-module structure, the following computation
shows that $\xi$ is $H$-equivariant. Let $x \in H$. Then 
\[
\begin{split}
\xi(m\cdot x)
& =\iota(m\cdot x) \otimes f = \rho^{-1}(x)  (\iota(m)\cdot x) \otimes
f  =
(\iota(m)\cdot x) \otimes \rho^{-1}(x)f \\
&
= ( \iota(m)\cdot x) \otimes (f\cdot x) =  \xi(m)\cdot x.
\end{split}
\]

We now prove that (4) implies (3). Let $\gamma \in
\operatorname{\mathcal X}(H)$ be a rational character of $H$ and
consider the rational
$H$-module $V = \Bbbk \gamma$. If $\xi$ and $W$ are as in condition  (4) we
conclude, using Lemma \ref{lema:1}, that $\gamma$ is extendible.

Assume now that (1) holds. Then by Theorem \ref{thm:istabcond}, $H$ is
observable in
$G$, and it follows from (2) that $\frac{1}{\det}\in E_M(H)$. Hence,
(5) holds.

Finally, if (5) holds, then by Theorem \ref{thm:obsercond},  $\mathcal
X(H) =E_G(H)$. Let
$\det\in \mathcal X (M)$ be a  determinant such that
$\frac{1}{\det}|_{_H}\in E_M(H)$, and let  $f\in \Bbbk[M]$ be an
extension of $\frac{1}{\det}|_{_H}$. If
$\chi\in \mathcal X(H)$, let $g\in \Bbbk[G]$ be an extension of
$\chi$. Then there exists $l\in \Bbbk[M]$ and $n\geq 0$ such that
$g=\frac{l}{\det^n}$. Therefore, for every
$x\in H$ and $a\in G$, we have that
\[
\chi(x)\frac{l(a)}{\det^n(a)}=\chi(x)g(a)=(g\cdot x)(a)=\frac{(l\cdot
  x)(a) }{\det^n(x)\det^n (a)}.
\]

It follows that $l\cdot a=\chi(x)\det^n(x) l$; that is $l\in \Bbbk[M]$
is an
extension of $\chi\det^n\in \mathcal X(H)$. Then $lf\in \Bbbk[M]$ is
an extension
of $\chi$.
 \qed \end{pf}

\begin{theorem}
\label{thm:obserrep}
Let  $M$ be an irreducible affine algebraic monoid with unit group $G$ and
let $H\subset G$ be a subgroup, closed in $M$. Then the following conditions
are equivalent.

\noindent (1) The subgroup $H$ is observable in $M$.

\noindent (2) There exists a finite dimensional rational right $M$-module
$V$ and $v\in V$ such that $G_v=H$.

In particular, if condition (2) holds, then $G_v=M_v$ is closed in $M$.
\end{theorem}

\begin{pf} 
 To prove that  condition  (1) implies condition (2), we
 adapt the proof of the case $M=G$ (see
 \cite[Corollary 7.3.6]{fer-ritt}). First we observe that, since $G$ is
 observable in $M$, $H$ is
 observable in $G$ and $\bigl[\Bbbk [G]^H\bigr]= \bigl[\Bbbk
 [G]\bigr]^H= \bigl[\Bbbk [M]\bigr]^H=\bigl[\Bbbk [M]^H\bigr]$.
Let $\{u_0, u_{1},\dots,u_{n}\} \subset \Bbbk[M]^H$ be such that
$\bigl\{\frac{u_1}{u_0},\dots,\frac{u_n}{u_0}  \bigr\}$ generates
 $ \bigl[\Bbbk [M]\bigr]^H$ over $\Bbbk$.  Let $W\subset \Bbbk[M]$ be the
finite dimensional rational right $M$-submodule generated by
$u_{0}, \dots, u_{n}$. Let
$V=\bigoplus_{i=0}^{n}W$,
$v_{0}=(u_{0}, \dots,u_{n}) \in V$
and consider the stabilizer $G_{v_{0}}$. It is clear that $H
 \subset G_{v_{0}}$.

 Conversely, if $a\in G_{v_0}$, then $u_{i}\cdot a =u_{i}$,
 $i=0,\dots,n$, and hence $\frac{u_{i}}{u_0}\cdot a=\frac{u_{i}}{u_0}$,
 $i=1,\dots,n$. As the elements
$\frac{u_{i}}{u_0}$,
$i=1,\dots,n$, generate  $\bigl[\Bbbk [M]\bigr]^H$, we
conclude that $f\cdot a=f$
for all $f\in \bigl[\Bbbk [M]\bigr]^H=\bigl[\Bbbk [G]\bigr]^H$. Hence,
$a\in H$, as follows by example from \cite[Corollary 7.3.4]{fer-ritt}.

To prove that (2) implies (1), let $v, V$ be as in (2) and
consider the following commutative diagram
\begin{center}
\mbox{
\xymatrix{
M\ar@{->}[r] &\overline{z\cdot M}\\
G\ar@{^(->}[u]\ar@{->}[r]& v\cdot G\ar@{^(->}[u]
}
}
\end{center}

Since $G/H$ is quasi-affine, it follows that $H$ is observable in $G$,
and hence $\bigl[\Bbbk[G]^H\bigr]=\bigl[ \Bbbk[G]\bigr]^H$. Thus,
\[
\bigl[\Bbbk[\overline{v\cdot M}]\bigr]=\Bbbk(v\cdot G)\cong
  \bigl[\Bbbk[G]^H\bigr]= \bigl[ \Bbbk[G]\bigr]^H.
\]

On the other hand, the orbit morphism $M\to \overline{v\cdot M}$ is
dominant, and hence
induces an inclusion $\iota: \Bbbk[\overline{v\cdot M}]\hookrightarrow
\Bbbk[M]$.   Since $G_v=H$, it follows that $f\bigl(v\cdot
(hm)\bigr)=f(v\cdot m)$ and hence $\iota\bigl(
\Bbbk[\overline{v\cdot M}]\bigr) \subset \Bbbk[M]^H$. Thus,
\[
\bigl[ \Bbbk[M]\bigr]^H = \bigl[ \Bbbk[G]\bigr]^H=
\bigl[\Bbbk[\overline{v\cdot M}]\bigr]\subset
\bigl[\Bbbk[M]^H\bigr]\subset   \bigl[ \Bbbk[M]\bigr]^H,
\]
and hence the action by right multiplication of $H$ on $M$ is observable.

Finally, recall that observable subgroups are closed in $M$.
 \qed \end{pf}

\begin{remark}[Open question]
\label{rem:Q1}
It is clear that all the results of this section remain valid when
considering right observability simply by adapting the statements and proofs
in an obvious way.

If $M=G(M)$ is an algebraic group, it is well known that $H\subset G$
is (left) observable in $G$ if and only if $H$ is right
observable. Indeed, in this case the antipode $S:\Bbbk[G]\to
\Bbbk[G]$, $S(f)(x)=f(x^{-1})$, induces an isomorphism $\Bbbk[G]^H\cong
{}^H\Bbbk[G]$. In the more general case of an algebraic monoid $M$, this line of
reasoning cannot be applied, since $S\bigr(\Bbbk[M]\bigl)$ is not
included in $\Bbbk[M]$. This raises the following question.

\medskip

\noindent  {\bf Q1} {\em Let $M$ be an algebraic monoid with unit group
$G$, and $H\subset G$ a closed subgroup, (left) observable in $M$. Is $H$ right
observable in $M$?}
\end{remark}

\section{Quotients of monoids by normal subgroups}
\label{sec4}

Let $M$ be an algebraic monoid of unit group $G$ and let $H\subset G$ a
closed subgroup, observable in $M$. Our goal is to prove that if
$H$ is normal in $G$, then the affinized quotient $M\aq H$ is
an affine monoid with unit group $G/H$. See Theorem
\ref{thm:normalimplisobse}. Before
doing so, we present a general result about the affinized quotient of
a monoid by an observable subgroup.

\begin{proposition}
Let $M$ be an algebraic monoid with unit group $G$ and let $H\subset G$ be a
closed subgroup, observable in $M$. Assume that $\Bbbk[M]^H$ is
finitely generated. Then $M\aq H$ is an affine embedding of $G/H$.
That is, $M\aq H$ is an affine $G$-variety, with an open orbit
isomorphic to $G/H$.
\end{proposition}
\begin{pf} 
Consider the following commutative diagram
\begin{center}
\mbox{
\xymatrix{
G\ar@{^(->}[r]\ar@{->>}[d]&M\ar@{->}[d]^-{\pi}\\
G/H\ar@{->}[r]_-{\varphi}& M\aq H
}
}
\end{center}
where the existence of $\varphi$ follows from the universal property of the
quotient. Since $H$ is observable in $M$, it follows from Theorem
\ref{thm:istabcond} that $\Bbbk(M\aq H)=\bigl[\Bbbk [M]^H\bigr]=
\bigl[ \Bbbk [M]\bigr]^H=\bigl[ \Bbbk [G]\bigr]^H=
\bigl[\Bbbk [G]^H\bigr]$. Hence, $\varphi$ is a birational dominant
$G$-morphism.  Since $G/H$ is a homogeneous $G$-space, it follows that
$\varphi$ is an open immersion.
 \qed \end{pf} 

\begin{proposition}
\label{prop:quotbynormred}
Let $M$ be a reductive monoid with unit group $G$, and let $H\subset G$ be a
normal subgroup, closed in $M$. Then
the categorical quotient $M\gq H$ is an affine algebraic monoid
with unit group $G/H$. Furthermore, $\pi: M\to M\gq H$ is a morphism of 
algebraic monoids.
\end{proposition}

\begin{pf} 
First of all, we prove that $H$ is observable in $M$. Indeed, $H$ is a
normal subgroup of $G$ and thus it
is a reductive group. Since $H\subset M$ is closed, it follows that
$gH$ is closed in $M$ for all $g\in G$, and we deduce from Proposition
\ref{prop:obsred} that  $H$ is observable in
$M$.

Let $(G\times  G)\times \Bbbk[M]\to \Bbbk[M]$, $\bigl((a,b)\cdot
f\bigr)(m)=f(a^{-1}mb)$ $a,b\in G$, $m\in M$, be the canonical action. If
$f\in \Bbbk[M]^H$ and  $c\in H$, then
\[
\begin{split}
\bigl( c\cdot \bigl((a,b)\cdot f\bigr)\bigr)(m) = & \ \bigl((c,1)\cdot
\bigl((a,b)\cdot f\bigr)\bigr)(m)=
f(a^{-1}c^{-1}mb)=\\
& \ f(la^{-1}mb)=f(a^{-1}mb)= \bigl((a,b)\cdot f\bigr)(m)
\end{split}
\]
where $l\in H$ is such that  $a^{-1}c^{-1}=la^{-1}$. If follows that
$\Bbbk[M]^H$ is a $(G\times G)$-submodule of $\Bbbk[M]$. Let now
$(c,d)\in H\times H$. Then $(1,d)\cdot f=f$. Indeed, if
$g\in G$, then $gd=lg$ for some $l\in H$, and hence
$f(gd)=f(lg)=f(g)$. In other words, $(1,d)\cdot f |_{_G}=f|_{_G}$, and thus
$(1,d)\cdot f =f$. It follows that $(c,d)\cdot f=f$. Hence, $M\gq H$
is an affine $(G/H\times G/H)$-variety and  the coproduct
$m^*:\Bbbk[M]\to \Bbbk[M]\otimes \Bbbk[M]$ is such that
\[
m^*\bigl(\Bbbk[M]^H\bigr)  \subset \bigl(\Bbbk[M]^H\otimes
\Bbbk[G]\bigr)\cap  \bigl(\Bbbk[G]
\otimes\Bbbk[M]^H\bigr)=\Bbbk[M]^H\otimes \Bbbk[M]^H.
\]
Moreover, since $\Bbbk[M]^H\subset \Bbbk[G]^H=\Bbbk[G/H]$, it follows
that  $M\gq H$ is an algebraic monoid and that we have a commutative
diagram
\begin{center}
\mbox{
\xymatrix{
G\ar@{^(->}[r]\ar@{->>}[d]&M\ar@{->>}[d]^-{\pi}\\
G/H\ar@{->}[r]_-{\varphi}& M\gq H
}
}
\end{center}
where the existence of $\varphi$ follows from the universal property of the
quotient. Then $\varphi$ is a dominant morphism of algebraic monoids, and
hence $G(M\gq H)=\varphi(G/H)$.
Since $H$ is observable in $M$ it
follows that
\[
\Bbbk(G/H)=\bigl[\Bbbk[G]\bigr]^H=\bigl[\Bbbk[M]\bigr]^H=
\bigl[\Bbbk[M]^H\bigr]   =\Bbbk(M//H).
\]
Thus,  $\varphi:G/H\to M\gq H$ is an injective  birational $G$-morphism,
and it  follows  that $\varphi:G/H\to G( M\gq H)$ is an open
immersion. Hence, $\varphi$ is an isomorphism.
 \qed \end{pf} 

The following result (with the exception of the last assertion)
follows from \cite[Theorem 2.5]{Re85} and \cite[Theorem 6.1]{Re05}.

\begin{theorem}
\label{thm:passtored}
Let $M$ be an affine algebraic monoid with unit group $G$. Then
$\Bbbk[M]^{R_u(G)}$ is a finitely generated algebra. Moreover,
$M\aq R_u(G)=\Spec\bigl( \Bbbk[M]^{R_u(G)}\bigr)$ is an affine
algebraic monoid with unit group $G/R_u(G)$. 
The morphism $\pi:M\to M\aq R_u(G)$ is a surjective morphism of
algebraic monoids, and satisfies the following universal property.
For any morphism $f:M\to N$, of algebraic monoids,
such that $f\bigl(R_u(G)\bigr)=\{1_N\}$ there exists a unique morphism 
$\widetilde{f}:M\aq R_u(G) \to N$ such that
the diagram
\begin{center}
\mbox{
\xymatrix{
M\ar@{->}[r]^-f\ar@{->>}[d]_-{\pi}&N\\
M\aq R_u(G)\ar@{->}[ru]_-{\widetilde{f}}&
}
}
\end{center}
commutes.
Moreover, the subgroup $H\subset G$ is closed in $M$ if and only if
$\pi(H)=\bigl(HR_u(G)\bigr)/R_u(G)\subset G/R_u(G)$ is
closed in $M\aq {R_u(G)}$.
\end{theorem}
\begin{pf} 
We provide a proof of the last assertion, since it
is not proved in  \cite{Re85} or \cite{Re05}.

It follows from \cite[Theorem 3.18]{Pu88} that $H=\overline{H}$ if and
only if $E(\overline{H})=\{1\}$.   But by \cite[Corollary 6.10]{Pu88},
$E(\overline{H})=\bigcup_{a\in H}aE(\overline{S})a^{-1}$, where
$S\subset H$ is a maximal tours of $H$. Hence, $H=\overline{H}\subset M$
if and only if $S=\overline{S}\subset M$ for any (some) maximal torus
$S\subset H$.

By \cite[Theorem 2.5]{Re85}, it follows that, if $T\subset G$ is a
maximal torus of $G$, then $\pi|_{_{\overline{T}}}: \overline{T}\to
\overline{\pi(T)}\subset M\aq R_u(G)$ is an isomorphism. Let now
  $S\subset H$ be a maximal torus of $H$  and $T\subset G$ a maximal torus
  of $G$ such that $S\subset T$. Then $\pi|_{_{\overline{S}}}:
  \overline{S}\to\overline{\pi(S)}\subset \overline{\pi(T)}$ is an
  isomorphism. In particular,
  $S=\overline{S}\subset M$ if and only if
  $\pi(S)=\overline{\pi(S)}\subset M\aq R_u(G)$. Since $\pi(S)\subset
  \pi(H)$ is a maximal tours, it follows that
  $H=\overline{H}\subset M$    if and only if
  $\pi(H)=\overline{\pi(H)}\subset M\aq R_u(G)$.

We conclude the proof by showing that the affinized quotient $\pi: M\to
M\aq R_u(G)$ 
is a surjective
morphism of algebraic monoids (this fact is implicit in
\cite{Re85,Re05} but was not stated or proved). By construction,
$\pi$ is a 
morphism of algebraic monoids. Since $M\aq R_u(G)$ is a reductive
monoid, it follows that  
\[
M\aq R_u(G)= G/R_u(G)E\bigl( M\aq
R_u(G)\bigr)G/R_u(G),
\]
see for example \cite[Theorem 4.2]{Re05}.
The surjectivity of $\pi$ follows now
from the fact that $\pi\bigl(E(M)\bigr)=E\bigl(M\aq R_u(G)\bigr)$.
 \qed \end{pf}

\begin{theorem}
\label{thm:normalimplisobse}
Let $M$ be an algebraic monoid with unit group $G$, and let $H\subset G$
be a normal subgroup, closed in $M$. Then the action $H\times M\to M$ is
observable. Moreover, if $\Bbbk[M]^H$ is finitely generated, then the affinized
quotient $M\aq H$ is an affine algebraic monoid,
with unit group $G/H$.
\end{theorem}
\begin{pf}  Let $\pi:M\to M\aq R_u(G)$ be the affinized quotient
as in Theorem \ref{thm:passtored}. By [ibid]
$M\aq R_u(G)$ is a reductive algebraic monoid with unit group
$G/R_u(G)$. We use $\pi$ to help us find a determinant function 
on $M$ (the function $\mu$ below)
with suitable properties. Consider the normal subgroup
$\pi(H)=\pi\bigl(HR_u(G)\bigr)\subset G/R_u(G)$. Then $\pi(H)$
is closed in $M\aq R_u(G)$ and hence by Proposition
\ref{prop:quotbynormred} it follows that $N= \bigl(M\aq
R_u(G)\bigr)\gq \pi(H)$ is an
affine algebraic monoid. Let $\rho: M\aq R_u(G)\to N$ be the
affinized quotient and $\chi:N\to \Bbbk$ a character such that
$\chi^{-1}(0)= N\setminus G(N)=\bigl(G/R_u(G)\bigr)/\pi(H)\cong
G/HR_u(G)$. Then $\mu: \chi\circ \rho\circ \pi:M\to \Bbbk$ is a character
such that $\mu(H)=1$ and $\mu^{-1}(0)= M\setminus G$. In particular
$\mu\in \Bbbk[M]^H$.   Since $H$ is normal in $G$ $H$
is observable in $G$. In other words, the action of $H$ on
$G=M_{\mu}$ is observable. It follows from Theorem
\ref{thm:obserchar}  that the action of $H$ on $M$ is observable.

To prove the last assertion we use the same arguments as in
 the proof of Proposition \ref{prop:quotbynormred}. Since
 $m^*\bigr(\Bbbk[M]^H\bigr) \subset \Bbbk[M]^H\otimes \Bbbk[M]^H$, it
 follows that
 $M\aq H$ is an
algebraic monoid and that we have a  commutative diagram
\begin{center}
\mbox{
\xymatrix{
G\ar@{^(->}[r]\ar@{->>}[d]&M\ar@{->}[d]^{\pi}\\
G/H\ar@{->}[r]_-{\varphi}& M\aq H
}
}
\end{center}
where the existence of $\varphi$ follows from the universal property of the
quotient. Then $\varphi$ is a dominant morphism of algebraic monoids, and
hence $G(M\aq H)=\varphi(G/H)$. Since the action of $H$ on $M$ is
observable, it follows from Theorem \ref{thm:istabcond} that
$\bigl[\Bbbk[M]^H\bigr]=\bigl[\Bbbk[M]\bigr]^H$, and hence $\varphi$ is a
birational morphism. Since $\rho$ is $G$-equivariant, to follows that
$\varphi$ is an open immersion,  and thus $G/H\cong G(M\aq H)$.
 \qed \end{pf}

If $H$ is normal in $G(M)$ then the affinized quotient $M\aq H$
satisfies a universal property in the category of algebraic monoids.

\begin{proposition} Let $N$ be an algebraic monoid and $f:M\to N$ a
  morphism of
  algebraic monoids such that $f(H)=\{1_N\}$. Then there exists a
  morphism of algebraic monoids  $\widetilde{f}:M\aq H\to N$ such that
  the following diagram commutes.

\begin{center}
\mbox{
\xymatrix{
M\ar@{->}[r]^-f\ar@{->>}[d]_-\pi&N\\
M\aq H\ar@{->}[ru]_-{\widetilde{f}}&
}
}
\end{center}
\end{proposition}

\begin{pf} 
Let
\begin{center}
\mbox{
\xymatrix{
1\ar@{->}[r]& L\ar@{->}[r]& N\ar@{->}[r]^-{\alpha}&A(N)\ar@{->}[r]& 0
}
}
\end{center}
be the Chevalley decomposition of $N$, where $\alpha$ is the Albanese
morphism  (see \cite[Theorem 3.2.1]{kn:brionlocal} and \cite[Theorem
1.1]{ritbr}). Since $M$ is affine, if follows that form the property
of the Albanese morphism that $f(M)\subset L$
(see  for example \cite[Theorem 5.1]{ritbr}). Since $f$ is a morphism
of algebraic monoids, such that $f(H)=\{1\}$, it follows that $f$ is
constant on the $H$-orbits, and hence $f^*:\Bbbk[L]\to \Bbbk[M]$ factors through
$\Bbbk[M]^H$. In other words, there exists $\widetilde{f}:M\aq H\to L$
such that $f=\widetilde{f}\circ\pi$.
 \qed \end{pf}

\begin{example}
\label{exam:normalimp}
Let
\[
H=\left\{ \left(\begin{smallmatrix} a^2 & 0\\
0& a^{-1}
\end{smallmatrix} \right) \mathrel{:} a\in \Bbbk^*\right\}\subset
\operatorname{GL}_2(\Bbbk) \subset \operatorname{M}_{2\times
  2}(\Bbbk)=M.
\]
Then $H$ is observable in $\operatorname{GL}_2(\Bbbk) =G(M)$ and
closed in $M$, but there exists no character $\chi:M\to \Bbbk$ such
that $H\subset \chi^{-1}(1)$ and $\chi^{-1}(0)= M\setminus G(M)$.
Since $H$ is closed in $G$ and since $\ell_g:M\to M$, $\ell_g(m)=g\cdot m$, is an
isomorphism, it follows that $gH\subset M$ is closed in $M$ for all $g\in
\operatorname{GL}_2(\Bbbk)$. Since  $H$ is reductive,
it follows from Proposition \ref{prop:obsred} that the action of $H$
on $M$ is observable.
\end{example}

\begin{remark}[Open question]
\label{rem:openquest}
Example \ref{exam:normalimp} shows that the condition of $H\subset G$
being a normal subgroup is crucial in Theorem
\ref{thm:normalimplisobse}. However, in that example the action of $H$ on
$M$ is observable. This raises the following question.

\medskip

\noindent {\bf Q2} {\em Let $M$ be an algebraic monoid with unit group
  $G$, and let $H\subset G$ an observable group, closed in $M$.
  Is the action $H\times M\to M$ observable?}
\end{remark}

\begin{remark}
If Q2 has a positive answer, then Q1 (see Remark \ref{rem:Q1}) has a
positive answer.
Indeed, assume that Q2 has a positive answer and let $H\subset M$ be
a left observable subgroup. Then by Theorem \ref{thm:istabcond} it follows
that $H$ is observable in $G$, and hence $H$ is right  observable in
$G$. If then follows from Q2 that $H$ is right observable in $M$.
\end{remark}


\begin{thebibliography}{100}



\bibitem{kn:obsdef}
 A.~Bialynicki-Birula, G.~Hochschild, G.D.~Mostow, \emph{Extensions of
  representations of algebraic linear groups}, Amer. J. Math. 85
(1963), 131-–144.

\bibitem{kn:brionlocal} M.~Brion, \emph{The local structure of algebraic
    monoids}, preprint, arXiv:0709.1255 [math.AG].

\bibitem{ritbr} M.~Brion, A.~Rittatore,
\emph{The structure of normal algebraic monoids},  Semigroup Forum  74
(2007),  no. 3, 410--422, arXiv:math/0610351 [math.AG].






\bibitem{fer-ritt} W.~Ferrer-Santos, A.~Rittatore, \emph{Actions and
    Invariants of Algebraic Groups.}
Series: Pure and Applied Mathematics,  268, Dekker-CRC Press,
Florida, (2005).





\bibitem{kn:grosshanslocal} F.~Grosshans, \emph{Localization and
    Invariant Theory}, Adv.~in Math. Vol.~21, No.~1 (1976) 50--60.






\bibitem{Pop} V.L.~Popov, {\em Stability criterion for the action
of a semisimple group on a factorial manifold}, Math. USSR Izv., 1970,
4 (3), 527--535.

\bibitem{Pu88} M.~S.~Putcha,
\emph{Linear Algebraic Monoids},
London Math. Soc. Lecture Notes Series {\bf 133},
Cambridge University Press, Cambridge, 1988.


\bibitem{Re85} L.~E.~Renner, \emph{Reductive monoids are von Neumann regular},
J. Algebra 93 (1985), no. 2, 237--245.


\bibitem{Re05} L.~E.~Renner,
\emph{Linear Algebraic Monoids},
Encyclop\ae dia of Mathematical Sciences {\bf 134},
Invariant Theory and Algebraic Transformation Groups, V,
Springer-Verlag, Berlin, 2005.

\bibitem{kn:oaoag} L.~E.~Renner, A.~Rittatore,
\emph{Observable actions of algebraic groups}, Preprint.
arXiv: 0902.0137v2 [math.AG].

\bibitem{Ri98} A.~Rittatore,
\emph{Algebraic monoids and group embeddings},
Transformation Groups {\bf 3}, No.~4 (1998), 375--396,
arXiv:math/9802073 [math.AG].

\bibitem{Ri06} A.~Rittatore,
\emph{Algebraic monoids with affine unit group are affine},
Trans. Groups, v. 12 3,  2009,  601--605, arXiv: math.AG/0602221. 







\end{thebibliography}
\end{document}